\documentclass[12pt,a4paper]{amsart}
\usepackage{amssymb,amsmath,amsthm}
\usepackage{graphicx}

\usepackage{epsfig}

\long\def\onefigure#1#2{%  #1 picture,  #2  caption
\begin{figure*}[tbp]
\begin{center}
#1
\end{center}
\caption{#2}
\end{figure*}
}%end onefigure def

\newcommand{\lipefig}[2]  % labeled Ipe figure
{\onefigure{\mbox{\psfig{file=#1.eps}}}{\label{f:#1} #2} }

\newtheorem{theorem}{Theorem}[section]
\newtheorem{lemma}{Lemma}[section]

\newcommand{\D}{\Delta}

\newcommand{\la}{\lambda}

\newcommand{\eps}{\varepsilon}

\newcommand{\R}{\mathbb{R}}
\newcommand{\N}{\mathbb{N}}
\newcommand{\Z}{{\mathbb{Z}^2}}
\newcommand{\A}{\mathcal{A}}
\newcommand{\T}{\mathcal{T}}
\newcommand{\F}{\mathcal{F}}

\newcommand{\p}{\mathcal{P}}

\newcommand{\remove}[1]{}

\numberwithin{equation}{section}

\begin{document}

\title{A matrix version of the Steinitz lemma}
\author{Imre B\'ar\'any}
\date{}							% Activate to display a given date or no date

\subjclass[2020]{Primary 05B20, secondary 52A22, 15A39}
\keywords{Vector sequences, their rearrangements, sign sequences.}

\maketitle

\begin{abstract}
The Steinitz lemma, a classic from 1913, states that $a_1,\ldots,a_n$, a sequence of vectors in $\R^d$ with $\sum_{i=1}^n a_i=0$, can be rearranged so that every partial sum of the rearranged sequence has norm at most $2d\max \|a_i\|$. In the matrix version $A$ is a $k\times n$ matrix with entries $a_i^j \in \R^d$ with $\sum_{j=1}^k\sum_{i=1}^na_i^j=0$. It is proved in \cite{OPW} that there is a rearrangement of row $j$ of $A$ (for every $j$) such that the sum of the entries in the first $m$ columns of the rearranged matrix has norm at most $40d^5\max \|a_i^j\|$ (for every $m$). We improve this bound to $(4d-2)\max \|a_i^j\|$.
\end{abstract}

\section{Introduction}

Let $\A_{d,k,n}$ denote the set of $k \times n$ matrices $A=\{a_i^j\}$ with $a_i^j \in \R^d$, and write $\sigma_m(A)$ for the sum of the entries in the first $m$ columns of $A$, that is, $\sigma_m(A)=\sum_{j=1}^k\sum_{i=1}^ma_i^j$. Assume $K$ is a $0$-symmetric convex body in $\R^d$, that defines, as usual, a norm $\|\cdot\|$ of $\R^d$. As $d,k,n$ are fixed in most cases we will simply write $\A$ for $\A_{d,k,n}$. In particular, $\A(K)$ stands for those $A\in \A$ satisfying $a_i^j \in K$ for every $i,j$, in other words $\|a_i^j\|\le 1$.

\medskip
Let $\{c_1^j,\ldots,c_n^j\}$ be a rearrangement of the vectors $\{a_1^j,\ldots,a_n^j\}$ in row $j$, for every $j\in [k]:=\{1,\ldots,k\}$ and write $C=\{c_i^j\}\in \A(K)$ for this row-wise rearranged, or row-permuted matrix. The target in the matrix version of the Steinitz lemma is to show that for a suitably permuted matrix $C$ all $\|\sigma_m(C)\|$ are bounded by a constant independent of $k$ and $n$. For this end $\sigma_n(A)=\sigma_n(C)$ has to be bounded, and we assume that $\sigma_n(A)=0$. Under this condition Oertel, Paat, and Weismantel \cite{OPW} showed the remarkable result that for every $A\in \A(K)$ there is a row-permuted copy, $C$, of $A$ with $\|\sigma_m(C)\|$ bounded, for every $m$, by a constant that only depends on $d$ and applied this result to certain proximity result for specially structured integer programs.  More precisely they proved the following theorem which they call the colourful Steinitz lemma.

\begin{theorem}\label{th:colour} For every $A \in \A(K)$ with $\sigma_n(A)=0$ there is a row-permuted copy, $C$, of $A$ such that
\[
\|\sigma_m(C)\| \leq \min \{dk, 40d^5\} \mbox{ for every } m\in [n].
\]
\end{theorem}

The $k=1$ case of this theorem is about vector sequences and is known as the Steinitz lemma~\cite{Stein} which states the following.

\begin{theorem}\label{th:stein} Assume $a_1,\ldots,a_N \in K$ and $\sum_{i=1}^N a_i=0$. Then there is a reordering $c_1,\ldots,c_N$ of the sequence $a_1,\ldots,a_N$ such that $\|\sum_{i=1}^mc_i \|\leq d$ for every $m \in [N]$.
\end{theorem}

The bound in Steinitz's 1913 paper \cite{Stein} was $2d$, the improved bound $d$ is an elegant result of Grinberg and Sevastyanov \cite{GS} which is more general because it holds for {\sl non-symmetric seminorms} as well. The definition is given in Section \ref{sec:nonsym} together with an example showing that Theorem \ref{th:colour} and Theorem~\ref{th:newver} below do not extend to this type of norms. 

\medskip
We remark that the case $k=1$ of Theorem~\ref{th:colour} is just the Steinitz lemma. The bound $dk$ in Theorem~\ref{th:colour} follows from the Steinitz lemma because the column sums of $A$ are vectors in $k$ times $K$ 
and their sum is zero; so this bound holds for non-symmetric seminorms as well.

\medskip
The main result of this paper is an improvement on the constant in Theorem~\ref{th:colour} with a different and simpler proof that gives a slightly more general result, namely Theorem~\ref{th:transf} below.

\begin{theorem}\label{th:newver} For every $A \in \A(K)$ with $\sigma_n(A)=0$ there is a row-permuted copy, $C$, of $A$ such that
\[
\| \sigma_m(C)\| \leq \min \{dk, 4d-2\} \mbox{ for every } m\in [n].
\]
\end{theorem}

It is interesting to note that the bounds in the above theorems are independent of the norm. The bound in the last theorem is at least $d/2$ for the $\ell_1$ norm implying that the upper bound in Theorem~\ref{th:newver} is of the right order of magnitude. The example showing this is given in  Section 4.

\section{Preparations}

We will need the following result which is due to B\'ar\'any and Grinberg~\cite{BG}.
\begin{theorem}\label{th:signs} If $a_1,\ldots,a_N \in K$, then there are signs $\varepsilon_1,\ldots,\varepsilon_N \in \{-1,1\}$ such that
$\|\sum_{i=1}^m \varepsilon_i a_i\| \le 2d-1$ for every $m \in [N]$.
\end{theorem}

\medskip
Suppose next that $B=\{b_i^j\} \in \A(K)$ and $\varepsilon=\{\varepsilon_i^j\}$ is a $k \times n$ matrix with $\varepsilon_i^j\in \{-1,1\}$ for all $j,i$. Writing $B^{\varepsilon}=\{\varepsilon_i^j b_i^j\} \in \A(K)$, define $V(B)=\min_{\varepsilon}\max_{m\in[n]}\|\sigma_m(B^{\varepsilon})\|$. This is the ``best sign assignment" to the matrix $B$. Set, finally
\[
V(K)=\max \{V(B): B\in \A(K)\}.
\]

\medskip
Theorem~\ref{th:signs} can be used to bound $V(K)$:

\begin{lemma}\label{l:S<SS} $V(K)\leq 2d-1$.
\end{lemma}

The {\bf proof} is simple. Consider the following ordering of the entries of $B$: $b_1^1,b_1^2\ldots,b_1^k,b_2^1,\ldots,b_2^k, \ldots,$ $b_{n-1}^k, b_n^1\ldots,b_n^k$ and apply Theorem~\ref{th:signs} to this sequence. You get signs $\varepsilon_i^j$ such that all partial sums with the $\varepsilon_i^j$ coefficients have norm at most $2d-1.$ Every $\sigma_m(B^{\varepsilon})$ is a partial sum of this form, so $V(B)\leq 2d-1.$\qed

\medskip
For $A \in \A(K)$ with $\sigma_n(A)=0$, set $U(A)=\min_C \max_{m\in [n]} \|\sigma_m(C)\|$ where the minimum is taken over all row-permuted copies, $C$, of $A$, and set
\[
U(K)= \max \{U(A): A\in \A(K) \mbox{ with }\sigma_n(A)=0\}.
\]
Taking the maximum instead of supremum is justified by compactness.

\medskip
Another ingredient of the proof is Chobanyan's transference theorem \cite{Chob}, see also \cite{Bar}. It is about sequences of vectors $a_1,\ldots,a_N\in K\subset \R^d$ so it is case $k=1$ of $U(K)$ and $V(K)$ and says that $U(K)\leq V(K)$. In other words, the ``best Steinitz constant'' of $K$ is not larger than the ``best sign assignment constant'' of $K$. We will not use this theorem directly, we use instead its proof scheme, with the necessary modifications. Our main result is another transference theorem, namely the following.

\begin{theorem}\label{th:transf} For every $k\ge 1$, $U(K)\leq 2V(K).$
\end{theorem}

This result together with Lemma~\ref{l:S<SS} clearly implies Theorem~\ref{th:newver}.

\section{Proof of Theorem \ref{th:transf}}

Let $A\in \A$ be a $k\times n$ matrix achieving the maximum in the definition of $U(K)$ and assume that the maximum is reached on the ordering $a_1^j,a_2^j,\ldots,a_n^j$ for every $j\in [k].$ This means that $U(A)=\max \{\|\sigma_m(A)\|: m\in [n]\}$ and for every row-permuted copy $C$ of $A$, $U(K)\leq \max \{\|\sigma_m(C)\|: m\in [n]\}$.

\medskip
Assume that $n$ is even, $n=2p$ say, the odd case will be treated later. We introduce an auxiliary matrix $B=\{b_i^j\} \in \A_{d,k,p}$ by setting $b_i^j=\frac 12(a_{2i-1}^j-a_{2i}^j)$ for all $j \in [k]$ and $i\in [p].$ Then $B\in \A(K)$ and by the definition of $V(K)$ there is a $k \times p$ matrix $\varepsilon=\{\varepsilon_i^j\}$ with $\varepsilon_i^j \in \{-1,1\}$ such that
$\max_{m\in [n]}\|\sigma_m(B^{\varepsilon})\| \leq V(K)$.

\medskip
Using the signs $\varepsilon_i^j$ we define a new ordering of row $j$ of $A$ (for every $j\in [k]$). This ordering is going to be $c_1^j,\ldots,c_n^j$ resulting in the row-permuted copy, $C$, of $A$, of course $C\in \A$. The rule is that $a_{2i-1}^j$ and $a_{2i}^j$ go to positions $i$ and $n+1-i$: if $\varepsilon_i^j=1$, then $c_i^j=a_{2i-1}^j$ and $c_{n+1-i}^j=a_{2i}^j$ and if $\varepsilon_i^j=-1$, then $c_i^j=a_{2i}^j$ and $c_{n+1-i}^j=a_{2i-1}^j.$ This is essentially the ordering in Chobanyan's transference theorem. For $i\leq p$ the $c_i^j$ can be expressed as
\[
c_i^j=\frac 12 \left [(a_{2i-1}^j+a_{2i}^j)+\varepsilon_i^j(a_{2i-1}^j-a_{2i}^j)\right].
\]

The other elements in row $j$ of the matrix $C$ are of the form $c_{n+1-i}^j$ where $i\le p$ and
\[
c_{n+1-i}^j=\frac 12 \left [(a_{2i-1}^j+a_{2i}^j)-\varepsilon_i^j(a_{2i-1}^j-a_{2i}^j)\right].
\]

\medskip
For $m\leq p$ we have
\begin{eqnarray*}
\sum_{i=1}^m c_i^j&=&\frac 12\left[ \sum_{i=1}^{m}(a_{2i-1}^j+a_{2i}^j) +\sum_{i=1}^m \varepsilon_i^j(a_{2i-1}^j-a_{2i}^j)\right]\\
  &=&\frac 12\sum_{i=1}^{2m}a_i^j+\sum_{i=1}^m \varepsilon_i^jb_i^j.
\end{eqnarray*}
Then, still for $m \leq p$,
\[
\sigma_m(C)= \frac 12\sum_{j=1}^k \sum_{i=1}^{2m}a_i^j+\sum_{j=1}^k \sum_{i=1}^m \varepsilon_i^jb_i^j
\in  \left[\frac 12 U(K)+V(K)\right]K
\]
because $\sum_{j=1}^k \sum_{i=1}^{2m}a_i=\sigma_{2m}(A)$ and $\sum_{j=1}^k\sum_{i=1}^m \varepsilon_i^jb_i^j=\sigma_m(B^{\varepsilon})$ and because $\|\sigma_{2m}(A)\|\leq U(A)=U(K)$ and $\|\sigma_m(B^{\varepsilon})\|\leq V(B)\leq V(K)$. This implies that $\|\sigma_m(C)\|\leq  \frac 12 U(K)+V(K)$.

\medskip
Now comes the case $m>p$. Since $\sigma_n(A)=\sigma_n(C)=0$ and $\sigma_n(C)=\sum_{j=1}^k\sum_{i=1}^nc_i^j=0$ we have
\begin{eqnarray*}
\sigma_m(C) &=&\sum_{j=1}^k\sum_{i=1}^mc_i^j=-\sum_{j=1}^k\sum_{i=m+1}^nc_i^j=-\sum_{j=1}^k\sum_{i=1}^{n-m}c_{n+1-i}^j\\
&=&-\sum_{j=1}^k\sum_{i=1}^{n-m}\frac 12 \left [(a_{2i-1}^j+a_{2i}^j)-\varepsilon_i^j(a_{2i-1}^j-a_{2i}^j)\right]\\
&=&-\sum_{j=1}^k \left[\frac 12\sum_{i=1}^{2(n-m)}a_i^j-\sum_{i=1}^{n-m}\varepsilon_i^jb_i^j \right]\\
&=&-\frac 12 \sum_{j=1}^k \sum_{i=1}^{2(n-m)}a_i^j+\sum_{j=1}^k\sum_{i=1}^{n-m}\varepsilon_i^jb_i^j,
\end{eqnarray*}
and $\|\sigma_m(C)\|\leq \frac 12 U(K)+V(K)$ follows the same way as above.

\medskip
As $C$ is a row-permuted copy of $A$, $\|\sigma_m(C)\|\geq U(A)$ for some $m \in [n]$ implying that
$U(K)\leq \|\sigma_m(C)\|\leq \frac 12 U(K)+V(K).$ This shows, in turn, that $U(K) \leq 2V(K)$.

\medskip
The odd case $n=2p+1$ is similar. For $i\leq p$ we define the same way $b_i^j=\frac 12(a_{2i-1}^j-a_{2i}^j)$ and set $b_{p+1}^j=a_n^j$. The matrix $B=\{b_i^j\} \in \A_{k,p+1}$ is again in $\A(K)$. We get signs $\varepsilon_i^j$ the same way as in the even case. In the new row-permuted copy, $C$, of $A$
the entries $a_{2i-1}^j$ and $a_{2i}^j$ go again to positions $i$ and $n+1-i$ just like before, and $c_{p+1}^j$ is going to be $\pm a_n^j$, does not matter which one. The partial sums $\sum_{i=1}^m c_i^j$ are estimated analogously. For $m\leq p$ there is no change, while for $m>p$,  $\sum_{j=1}^k\sum_{i=1}^mc_i^j=-\sum_{j=1}^k\sum_{i=m+1}^n c_i^j$ and $a_n^j$ does not appear in any one of these partial sums. Again $\sigma_m(C)\in  \left[\frac 12 U(K)+V(K)\right]K$ for all $m \in [n]$ and the proof is finished the same way as in the even case.\qed

\section{The example giving a lower bound on $U(K)$} 

The example is for the $\ell_1$ norm and gives a matrix $A$ with $U(A)\geq \frac d2.$

\medskip
We start with a $d \times 2d$ matrix $A\in \A_{d,d,2d}$ whose rows are identical and are of the form
$e_1,\ldots,e_d, -\frac ed,\ldots,-\frac ed$ where $e_1,\ldots,e_d$ are the standard basis of $\R^d$ and $e=e_1+\ldots+e_d$. It is clear that $\sigma_{2d}(A)=0$

\medskip
Assume that in a row-permuted copy $C$ of $A$ a fixed column contains $q$ copies of $-\frac ed$, the other entries are $d-q$ basis vectors. The $\ell_1$ norm of the sum of the entries in this column is at least
\[
(d-q)\left(1-\frac qd\right)+ q\frac qd=\frac {(d-q)^2+q^2}d \geq \frac d2,
\]
as one can check directly. So the sum of the entries in every column of $C$ has norm at least $\frac d2$ and then $\|\sigma_1(C)\| \geq \frac d2$. Consequently $U(A)\geq \frac d2$.

\medskip
We take $s$ copies of this matrix to get a new matrix in $\A_{d,d,2sd}$. Adding $k-d$ all zero rows to it, the resulting matrix $B\in\A_{d,k,n}$ with $n=2ds$ satisfies $\sigma_n(B)=0$. The previous argument shows that the sum of the entries in every column of a permuted copy of $B$ have $\ell_1$ norm at least $\frac d2$,
implying that $U(B)\geq \frac d2.$

\section{Non-symmetric seminorms}\label{sec:nonsym}

A (non-symmetric) seminorm is given by a compact convex set $K \subset \R^d$ containing the origin in its interior, and the norm of $x\in \R^d$ is given by $\|x\|=\|x\|_K=\min \{t\ge 0: x\in t \cdot K\}.$ Such a norm satisfies the triangle inequality $\|x+y\|\le \|x\|+\|y\|$ and the equality $\|\alpha x\|= \alpha \|x\|$ for all $\alpha \ge 0$. Recall that for usual norms the equality $\|\alpha x\|= |\alpha| \|x\|$ holds for all $\alpha \in \R$. 

\medskip
For the example showing that Theorem~\ref{th:colour} and Theorem~\ref{th:newver} do not extend to seminorms we define $K$ as the set of points $x=(x_1,\ldots,x_d)\in \R^d$ satisfying $-1\le x_i$ for all $i \in [d]$ and $\sum_{i=1}^dx_i\le M$ where $M>0$ is large,  for instance $M=2^d$ would do.

\medskip
The matrix $A$ in the example is going to be a $k \times d$ matrix where $k$ is approximately $2^d/d$. The entries $a_i^j\in \R^d$ for all $i,j \in [d]$ are defined by $a_i^j=2^{j-1}e_i-e$ and form a $d\times d$ matrix, to be denoted by  $A^*$. For the rest of the entries of $A$ the preliminary value of $a_i^j$ is going to be $b_i^j=-e$ that will be modified slightly soon. Direct checking shows that 
\[
\sum_{j=1}^d\sum_{i=1}^da_i^j=(2^d-1-d^2)e,\mbox{ and }\sum_{j=d+1}^k\sum_{i=1}^db_i^j=-(k-d)de.
\]
As $\sigma_d(A)$ has to be zero we have to change some of the $b_i^j$s. 

Define first $k=\lceil\frac{2^d-1}d\rceil$ and $s$ via $2^d-1=kd-s$ so $0\le s <d$.  Set next $f=e_1+\ldots+e_{\lfloor d/2\rfloor}$ and $g=e_{\lfloor d/2\rfloor+1}+\ldots+e_d$ and define $a_i^{d+1}=-f$ and $a_i^{d+2}=-g$ for $i\in [s]$. For the rest of the entries we set $a_i^j=b_i^j=-e$. Then $\sigma_d(A)=0$ and $\|a_i^j\|\le 1$ holds for every $i,j$ (with equality in fact).

\medskip
Let $C$ be a row-permuted copy of $A$.  Write $S$ for the sum of the entries in column one of $C$.  Assume first that this column does not contain any one of entries of some column of $A^*$, of column $h$, say. Then the coefficient of $e_h$ in $S$ is $-k$ or $-k+1$ depending on what entry column one of $C$ contains from rows $d+1$ and $d+2$.  Consequently $\|S\|\ge k-1$.

Assume next that column one of $C$ contains at least one entry from every column of $A^*$.  Then it has to contain exactly one entry from every column of $A^*$. Let $c_1^1$ come from column $h$ of $A^*$, so $c_1^1=a_h^1$. Then the coefficient of $e_h$ in $S$ is $1-k$ or $1-k+1$ depending again on what entry column one of $C$ contains from rows $d+1$ and $d+2$.  Consequently $\|S\| \ge k-2$.

Thus $\|S\| \ge k-2$ in both cases and since $\sigma_1(C)=S$ we have 
\[
U(A)\ge \|\sigma_1(C)\| \ge \lceil\frac{2^d-1}d\rceil-2\ge \frac {2^d}{2d}\mbox{ when } d\ge 4.
\]

\medskip
Together with the bound $U(K)\le dk$ this gives the estimate
\[ 
 k-2 \le U(K) \le dk,
\]
so $U(K)$ somewhere is between $\frac {2^d}{2d}$ and $2^d$.

\bigskip
{\bf Acknowledgements.} I am indebted to an anonymous referee and to several colleagues for careful reading an earlier version of this paper and for their comments. Special thanks are due to Marco Caoduro for an observation that reduced the previous bound on $U(K)$ by a factor of two. This piece of work was partially supported by Hungarian National Research Grants No 131529, 131696, and 133819.

\bigskip
\vspace{.5cm} {\sc Imre B\'ar\'any}

\medskip
\noindent
Alfr\'ed R\'enyi Institute of Mathematics\\
13-15 Re\'altanoda Street, Budapest, 1053 Hungary,\\
e-mail: {\tt imbarany@gmail.com}\\
and

\vskip.1cm
\noindent
Department of Mathematics, University College London\\
Gower Street, London WC1E 6BT, UK

\end{document}